\documentclass{article}

\usepackage{arxiv}

\usepackage[utf8]{inputenc} 
\usepackage[T1]{fontenc}    
\usepackage{hyperref}       
\usepackage{url}            
\usepackage{booktabs}       
\usepackage{amsfonts}       
\usepackage{nicefrac}       
\usepackage{microtype}      
\usepackage{lipsum}		
\usepackage{graphicx}
\usepackage{natbib}
\usepackage{doi}
\usepackage{caption}
\usepackage{csvsimple}
\usepackage{booktabs, tabularx, wrapfig}
\usepackage{array, ltablex, multirow}
\usepackage{placeins}
\usepackage{graphicx}
\usepackage{subcaption}
\usepackage{wrapfig,floatrow}
\usepackage[nonumberlist,nogroupskip]{glossaries}

\setglossarystyle{super}

\hypersetup{%
  colorlinks=false,
  linkbordercolor=[rgb]{0.6,0.6,0.6},
  pdfborderstyle={/S/U/W 1}
}

\title{ François Viète and
his contribution to mathematics}

\author{ \href{https://orcid.org/0000-0001-6845-8964}{\includegraphics[scale=0.06]{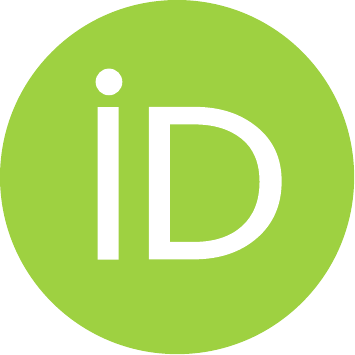}\hspace{1mm}Athanasios Paraskevopoulos}\thanks{Author earned a Bachelor's degree focused in Mathematics from the University of Aegean. Currently MSc. student in Mathematics at Hellenic Open University. He seeks to engage his students in a search for relating mathematics to real life's problem. } \\
	School of Science and Technology,\\Hellenic Open University, 263 35  Patras, Greece\\
	\texttt{std156374@ac.eap.gr} \\
}

\renewcommand{\shorttitle}{\textit{arXiv} Template}

\hypersetup{
pdftitle={A template for the arxiv style},
pdfsubject={q-bio.NC, q-bio.QM},
pdfauthor={David S.~Hippocampus, Elias D.~Striatum},
pdfkeywords={First keyword, Second keyword, More},
}

\begin{document}
\maketitle

\begin{abstract}
This paper studies the work of the French mathematician François Viète, known as the ”father of modern algebraic` notation”. Along with this fundamental change in algebra, Viète adopted a radically new notation based on Greek geometric` equalities. Its letters represent values rather than types, and its given values are undefined. Where algebra had previously relied on polynomials as sets, Viète became the first modern algebraist to work with polynomials generated by operations, and the notations reflect these notions. His work was essential to his successors because it enabled those mathematicians who followed him to develop the mathematics we use today.
\end{abstract}


\section{François Viète - Life}

\begin{wrapfigure}{r}{0.5\textwidth}
  \vspace{-20pt}
  \begin{center}
    \includegraphics[width=0.48\textwidth]{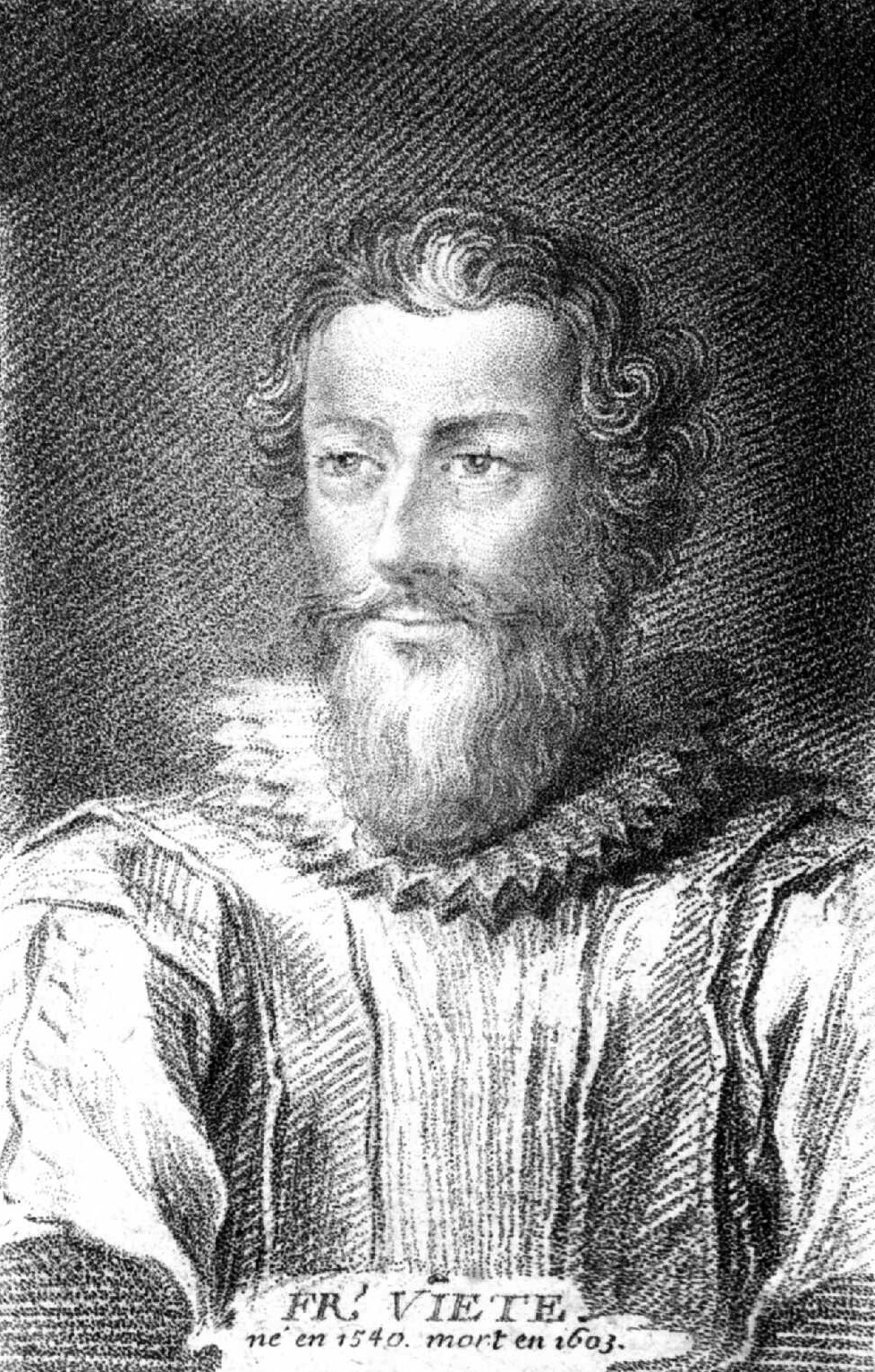}
  \end{center}
  \vspace{-2pt}
  \caption{Portrait of Viète from Savérien’s 1773 Histoire des Philosophes.}
  \vspace{-10pt}
\end{wrapfigure}

François Viète (or Vieta, as he is often known by his latinized name) was born in 1540 in Fontenay-le-Comte, France. His grandfather was a merchant from La Rochelle. His father, Etienne Viète, was an attorney in Fontenay-le-Comte and a notary in Le Busseau. His mother was the aunt of Barnabé Brisson, a magistrate and the first president of parliament during the ascendancy of the Catholic League of France. Educated in a a Franciscan school, then at the University of Poitiers,graduating as a Bachelor of Laws in 1559. He took up the practice of law in his home town.\cite{Cantor1911}

 He soon rose to prominence by his astute legal services to prominent people, and later served as royal councilor to Kings Henry III and Henry IV of France. In his spare time he worked on mathematics and published his results at his own expense. He has been called the father of modern algebra \cite{Brit2016} and the foremost mathematician of the sixteenth century \cite{Kline1972}.\\

A few stories will indicate something of his character. While working for King Henry III, he discovered the key to a Spanish cipher of 500 characters, and so was able to read the secret correspondence of his enemies. Philipp II of Spain was so sure that his code was invulnerable that when he heard of this, he complained to the Pope that the French were using sorcery against him, contrary to good Christian morals.\cite{Brit2016}\\

Viète's tact in dealing with people is illustrated by the case of Françoise de Rohan, a cousin of Henry III. She had been engaged to the duke J. de Nemours, and had a son by him, but then he married another woman, Anne d'Este. Françoise wanted him declared her legal husband, and his children by Anne declared bastards. Viète found a solution in that parliament declared her the legal wife of Nemours, and gave her a property as her dukedom, but they also declared the marriage dissolved, so that Anne d'Este and her children were not harmed either in honor or property.\cite{Vonwort1970}\\

 Viète's mathematical skill appears in the following incident from the summer of 1594. The Belgian mathematician A. van Roomen had proposed a challenge to all contemporary mathematicians, to solve a certain equation of the 45th degree. The Dutch ambassador presented van Roomen's book to King Henry IV with the comment that apparently there were no mathematicians of any importance in France. The King called for Viète who immediately found one solution to the equation, and then next day presented 22 more.\cite{Vonwort1970}\\
 
 In response to van Roomen, Viète challenged him to solve the problem of Apollonius, to construct a circle tangent to three given circles. When Adrianus Romanus wanted to effect a solution using two hyperbolas, Viète was little content with that solution, because it was alien to the usage of geometers, according to which in the solution of problems that are called plane, only circles and lines should be used. He himself gave a universal solution to that problem of tangencies by a more geometric method, and published a little book in 1600 in Paris, to which he gave the title Apollonius Gallus \cite{Witmer1983} 
 Adrianus was so pleased with this work that he immediately set out for France to see Viète, and struck up an intimate friendship with him. Viète received him hospitably, was together with him for a whole month, and took care of him to the extent of his resources.  \cite{Camerer1795}\\
 
 We may be surprised that such a busy lawyer as Viète should have had any time for mathematics. But he had a contemplative life as well. A contemporary historian said in 1620, his meditations were so profound that he was often seen fixed in cogitation for three days continuously, sitting at his eating table without food or sleep, except what he could get leaning on his elbow, and neither moved nor sought refreshment at natural intervals \cite{Shooten1646}\\
 
Viète died in 1603, less than three months after he had received the King's permission to retire from public life.

 \section{New Algebra}
 
At the end of the 16th century, mathematics was placed under the dual aegis of the Greeks, from whom it borrowed the tools of geometry, and the Arabs, who provided procedures for the resolution. At the time of Viète, algebra therefore oscillated between arithmetic, which gave the appearance of a list of rules, and geometry which seemed more rigorous. Meanwhile, Italian mathematicians Luca Pacioli, Scipione del Ferro, Niccolò Fontana Tartaglia, Ludovico Ferrari, and especially Raphael Bombelli (1560) all developed techniques for solving equations of the third degree, which heralded a new era.\\
\begin{figure}[ht]
\includegraphics[width=8cm]{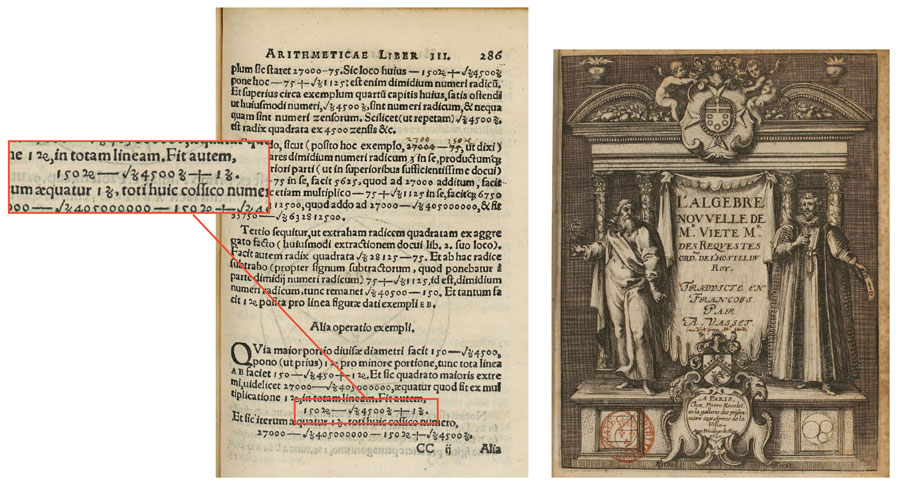}
\centering
\caption{Left: A polynomial from Michael Stifel’s book Arithmetica Integra (1544), showing premodern notation. We would write it as $150x-\sqrt(4500x^2)+x^2$. Note the coefficient of “1” on the last term. Compare with the notation on the preceding page. Right: Title page to Vasset’s 1630 French translation of two of Viète’s works, showing Viète on the right.}
\end{figure}

On the other hand, the German school of the Coss, the Welsh mathematician Robert Recorde (1550) and the Dutchman Simon Stevin (1581) brought an early algebraic notation, the use of decimals and exponents. However, complex numbers remained at best a philosophical way of thinking and Descartes, almost a century after their invention, used them as imaginary numbers. Only positive solutions were considered and using geometrical proof was common.\\

The task of the mathematicians was in fact twofold. It was necessary to produce algebra in a more geometrical way, i.e., to give it a rigorous foundation; and on the other hand, it was necessary to give geometry a more algebraic sense, allowing the analytical calculation in the plane. Viète and Descartes solved this dual task in a double revolution.

\subsection{The Foundations of Algebra} 
Viète gave algebra a foundation as strong as geometry. He then completed the algebra of processes (al-Jabr and Muqabala), creating the first symbolic algebra. In this way, he did not hesitate to say that with this new algebra, all problems could be solved (nullum non problema solvere). In his treatise In artem analyticam isagoge (Tours, 1591), Viète demonstrated the value of symbols by introducing letters to represent unknowns. He proposed the use of letters as symbols for quantities, both known and unknown. He used vowels for unknowns and consonants for known quantities. The convention where letters near the beginning of the alphabet represent known quantities while letters near the end represent unknown quantities was later introduced by Descartes in La Gèometrie.\cite{Cogito2013} This convention is used today, often without people realizing that a convention is used at all.\cite{Connor}
Also,  François Viète is credited with popularizing the use of the Maltese cross (+) as the arm symbol for addition and (–) for subtraction. This was first used by Johannes Widmann in a book in 1489, but only became popular when Viète started using them in his work. Previously, these were usually written as plus and minus in arithmetic operations. He also used the fraction line as a division symbol and the small word "in" as a fixed abbreviation of multiplication. Furthermore, Viète expressed the equality of two terms with the word "aequabitur" and thus invented the first equal sign. He also wrote related terms below each other and connected them with curly brackets.\\

\begin{figure}[ht]
\includegraphics[width=5cm]{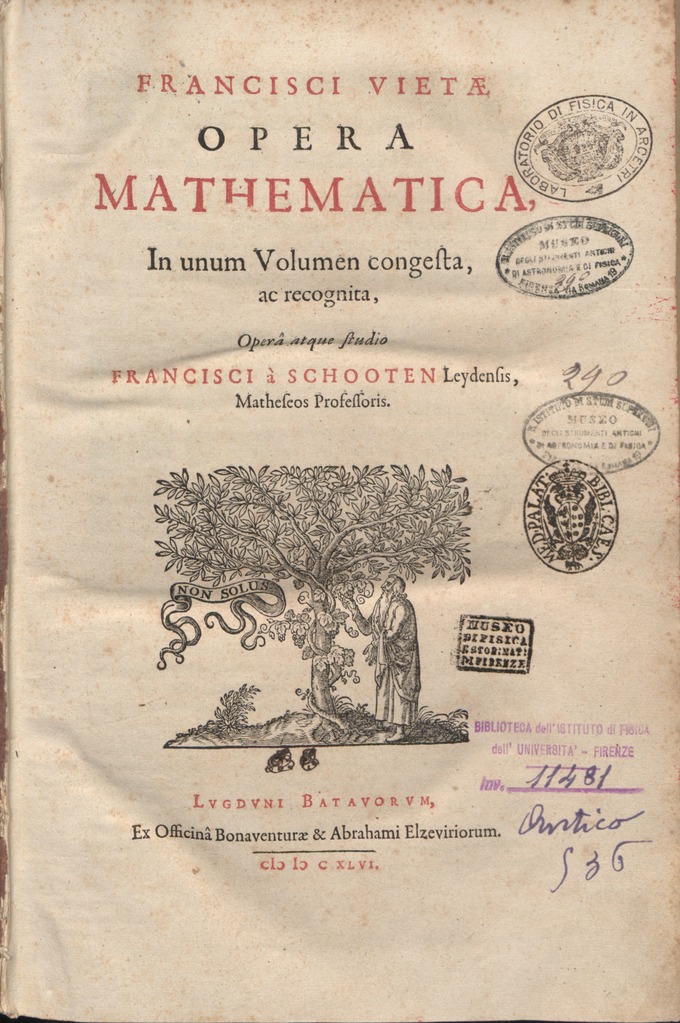}
\centering
\caption{Francisci Vietae Opera mathematica in unum volumen congesta ac recognita, opera et studio Francisci a Schooten, 1646. Leydensis.}
\end{figure}

 His contribution to the theory of equations is De aequationum recognitione et emendatione (1615; Concerning the Recognition and Emendation of Equations), in which he presented methods for solving equations of the second, third, and fourth degree. Viète knew the connection between the positive roots of an equation (which, in his time, were thought of as the only roots) and the coefficients of the different powers of the unknown quantity.\cite{Brit} Viéte’s new systematic algebra opened doors for further mathematical development enabling his successors to express their thoughts mathematically with a higher degree of detail and precision.\\
 
 Viète also remained a prisoner of his time in several respects. First, he was the heir of Ramus and did not address the lengths as numbers. His writing kept track of homogeneity, which did not simplify their reading. He failed to recognize the complex numbers of Bombelli and needed to double-check his algebraic answers through geometrical construction. Although he was fully aware that his new algebra was sufficient to give a solution, this concession tainted his reputation.\\

However, Viète created many innovations: the binomial formula, which would be taken by Pascal and Newton, and the coefficients of a polynomial to sums and products of its roots, called Viète's formula.

\subsection{Geometric algebra}
Viète was well skilled in most modern artifices, aiming at the simplification of equations by the substitution of new quantities having a certain connection with the primitive unknown quantities. Another of his works, Recensio canonica effectionum geometricarum, bears a modern stamp, being what was later called an algebraic geometry -- a collection of precepts on how to construct algebraic expressions with the use of a ruler and compass only. While these writings were generally intelligible, and therefore of the greatest didactic importance, the principle of homogeneity, first enunciated by Viète, was so far in advance of his times that most readers seem to have passed it over. That principle had been made use of by the Greek authors of the classic age; but of later mathematicians, only Hero, Diophantus, etc., ventured to regard lines and surfaces as mere numbers that could be joined to give a new number, their sum.\cite{Cantor1911}\\

The study of such sums, found in the works of Diophantus, may have prompted Viète to lay down the principle that quantities occurring in an equation ought to be homogeneous, all of their lines, or surfaces, or solids, or supersolids — an equation between mere numbers being inadmissible. During the centuries that have elapsed between Viète's day and the present, several changes of opinion have taken place on this subject. Modern mathematicians like to make homogeneous equations as are not so from the beginning, in order to get values of a symmetrical shape. Viète himself did not see that far; nevertheless, he indirectly suggested the thought. He also conceived methods for the general resolution of equations of the second, third, and fourth degrees different from those of Scipione dal Ferro and Lodovico Ferrari, with which he had not been acquainted. He devised an approximate numerical solution of equations of the second and third degrees, wherein Leonardo of Pisa must have preceded him, but by a method which was completely lost.\cite{Cantor1911}\\

Above all, Viète was the first mathematician who introduced notations for the problem (and not just for the unknowns)\cite{Bos}.As a result, his algebra was no longer limited to the statement of rules but relied on an efficient computational algebra, in which the operations act on the letters, and the results can be obtained at the end of the calculations by a simple replacement. This approach, which is the heart of the contemporary algebraic method, was a fundamental step in the development of mathematics\cite{Pycior}. With this, Viète marked the end of medieval algebra (from Al-Khwarizmi to Stevin) and opened the modern period.

\section{The logic of species}

Viète made significant contributions to trigonometry, algebra, and geometry.\\

His first published work, the Canon mathematicus \cite{Canon1579}  has trigonometric tables computed to 9 decimal places, and contains a systematic collection of trigonometric formulas. Because of many misprints and a misunderstanding with the editor, this volume was not included in his collected works.\\

The systematic presentation of his mathematic theory, which he called "species logistic" (from species: symbol) or art of calculation on symbols (1591).\\

He described in three stages how to proceed with solving a problem:
\begin{itemize}
    \item 

As a first step, he summarized the problem in the form of an equation.  called this stage the Zetetic. It denotes the known quantities by consonants (B, D, etc.) and the unknown quantities by the vowels (A, E, etc.)
\item In a second step, he made an analysis. He called this stage the Poristic. Here mathematicians must discuss the equation and solve it. It gives the characteristic of the problem, porisma, from which we can move to the next step.
\item In the last step, the exegetical analysis, he returned to the initial problem which presents a solution through a geometrical or numerical construction based on porisma.
\end{itemize}
Among the problems addressed by Viète with this method is the complete resolution of the quadratic equations of the form 
$X^{2}+Xb=c$ and third-degree equations of the form $X^{3}+aX=b$ (Viète reduced it to quadratic equations). \\

He knew the connection between the positive roots of an equation (which, in his day, were alone thought of as roots) and the coefficients of the different powers of the unknown quantity (see Viète's formulas and their application on quadratic equations). He discovered the formula for deriving the sine of multiple angles, knowing that of the simple angle with due regard to the periodicity of sines. This formula must have been known to Viète in 1593.\cite{Cantor1911}\\

Viète's most significant contributions were in algebra. While letters had been used to describe an unknown quantity by earlier writers, Viète was the first to also use letters for the parameters or constant coefficients in an equation. Thus, while Cardano solved particular cubic equations such as
$$x^3+6x=20$$

Viète could treat the general cubic equation\cite{Vonwort1970}

$$x^3+px=q$$

where $p$ and $q$ are constants.\\

Viète's algebra was significantly more systematic in the formal manipulation of equations than that of his predecessors, but it still does not reach the facility of modern techniques, because he did not consider negative numbers, and did not yet have a symbol for equality. For example, he writes the above cubic equation.

$$ \text{A cubus } +\text{ quad. in A, æquetur B quad. in Z.}$$
We would write this in symbols

$$A^3+B^2A=B^2Z.$$
He used vowels, such as $A$ here, to represent unknowns, and consonants, such as $B$ ,$Z$, to denote constants. Viète also preserved the old distinction between linear, planar, and solid quantities, so that $A^3$  represented the solid cube whose side is $A$. In equations, solids can only be compared with other solids, so what we call $p$,$q$ he had to call $B^2$  and $B^2Z$  to preserve homogeneity.\\

In geometry, Viète gave a solution to the problem of Apollonius referred to above, and also made a study of ``solid" problems such as the trisection of the angle and the construction of the regular heptagon, which use a marked ruler in addition to the Euclidean tools of ruler and compass.

\section*{Viète's formula  }

In algebra, Vieta's formulas are a set of results that relate the coefficients of a polynomial to its roots. In particular, it states that the elementary symmetric polynomials of its roots can be easily expressed as a ratio between two of the polynomial's coefficients.

\subsection*{Statement}
Let $P(x) = a_n x^n + a_{n-1}x^{n-1} + \cdots + a_1 x + a_0$ be any polynomial with complex coefficients with roots $r_1, r_2, \ldots , r_n$, and let $s_j$ be the $j^{\text{th}}$ elementary symmetric polynomial of the roots.\hfill \break

Vieta’s formulas then state that \[s_1 = r_1 + r_2 + \cdots + r_n = - \frac{a_{n-1}}{a_n}\]\[s_2 = r_1r_2 + r_1r_3 + \cdots + r_{n-1}r_n = \frac{a_{n-2}}{a_n}\]\[\vdots\]\[s_n = r_1r_2r_3 \cdots r_n = (-1)^n \frac{a_0}{a_n}.\] This can be compactly summarized as $s_j = (-1)^j \frac{a_{n-j}}{a_n}$ for some $j$ such that $1 \leq j \leq n$.
\subsection*{Proof}
Let all terms be defined as above. By the factor theorem, $P(x) = a_n (x-r_1)(x-r_2) \cdots (x-r_n)$. We will then prove Vieta’s formulas by expanding this polynomial and comparing the resulting coefficients with the original polynomial’s coefficients.\hfill \break

When expanding this factorization of $P(x)$, each term is generated by a series of $n$ choices of whether to include $x$ or the negative root $-r_{i}$ from every factor $(x-r_{i})$. Consider all the expanded terms of the polynomial with degree $n-j$; they are formed by multiplying a choice of $j$ negative roots, making the remaining $n-j$ choices in the product $x$, and finally multiplying by the constant $a_n$.\hfill \break

Note that adding together every multiplied choice of $j$ negative roots yields $(-1)^j s_j$. Thus, when we expand $P(x)$, the coefficient of $x_{n-j}$ is equal to $(-1)^j a_n s_j$. However, we defined the coefficient of $x^{n-j}$ to be $a_{n-j}$. Thus, $(-1)^j a_n s_j = a_{n-j}$, or $s_j = (-1)^j a_{n-j}/a_n$, which completes the proof.

\section*{Viète's formula for $\pi$}

In mathematics, Viète's formula is the following infinite product of nested radicals representing twice the reciprocal of the mathematical constant $\pi$:

$$\pi=2 \times \frac{2}{\sqrt{2}}\times\frac{2}{\sqrt{2+\sqrt{2}}}\times\frac{2}{\sqrt{2+\sqrt{2+\sqrt{2}}}}\times \cdots$$

 He calculated 10 decimal places of $\pi$ by applying the Archimedes method to a polygon with $6 \times 2^{16} = 393,216$ sides.

   \begin{figure}[ht]
\includegraphics[width=9cm]{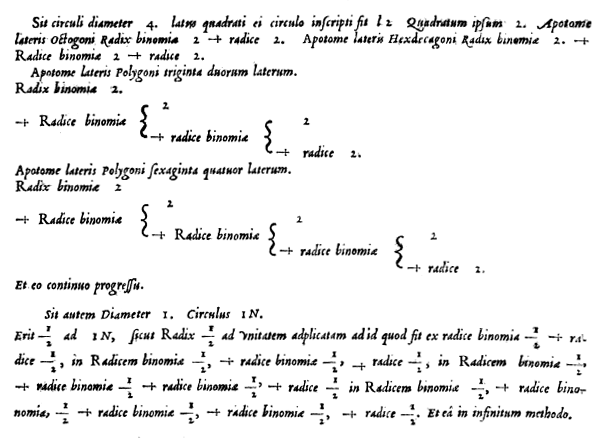}
\centering
\caption{Viète's formula, as printed in Viète's Variorum de rebus mathematicis responsorum, liber VIII (1593)}
\end{figure}
The formula is named after François Viète, who published it in 1593. As the first formula of European mathematics to represent an infinite process, it can be given a rigorous meaning as a limit expression and marks the beginning of the mathematical analysis. It has linear convergence and can be used for calculations of $\pi$, but other methods before and since have led to greater accuracy. It has also been used in calculations of the behavior of systems of springs and masses, and as a motivating example for the concept of statistical independence.\\

The formula can be derived as a telescoping product of either the areas or perimeters of nested polygons converging to a circle. Alternatively, repeated use of the half-angle formula from trigonometry leads to a generalized formula, discovered by Leonhard Euler, that has Viète's formula as a special case. Many similar formulas involving nested roots or infinite products are now known.

\section{Viète's construction of the regular heptagon }

Viète is best known for his work in algebra. However, I would like to describe here some of his work in geometry, which seems to be less well known: the trisection of the angle, the solution of the ``casus irreducibilis'' of the cubic equation, and the construction of the regular heptagon. These results are all presented in the Supplementum Geometriae, first published in 1593 (though I have not been able to locate any copy of the original), and included in Opera\cite{Shooten1646}.\hfill \break

Viète makes clear from the beginning that besides ruler and compass, he is admitting another operation:\hfill \break

In order to supplement the defect of geometry, let it be allowed from any point to any two straight lines, to draw a straight line cutting off between them any segment fixed in advance. \cite{Shooten1646}\hfill \break

In other words, given two lines l,m, given a point A, and a segment d, we can draw a line through A so that the two lines l,m cut off on it a segment BC=d. Viète also allows the same where a line and a circle are given instead of two lines.\hfill \break

One can imagine this construction using a marked ruler: one makes two marks on the ruler at a distance d apart, then slides the ruler so that the marks lie over the lines l,m while passing still through the point A. This process goes by the name of neusis in Greek.Viète does not give it a name.\hfill \break

Viète's work has 25 propositions, most of which I will skip over.\hfill \break

Proposition IX \cite{Shooten1646} is to trisect an angle. Let the given angle DBE be placed in a circle with B at the center. Draw a line through E, meeting BD extended at F, so that the segment FG cut off by the line and the circle is equal to BD (this is the neusis step). Then the angle at F is one third of the original angle at B. For the proof, draw BG. Call the angle at F 1. Then BFG is isosceles, so the angle GBF is also 1. The angle BGE is an exterior angle, hence equal to 2 times 1. BGE is isosceles, so GEB= 2. Now the original angle EBD is exterior to the triangle BEF so it is 3 times the angle at F.
\hfill \break

The question naturally arises, did Viète discover this result himself, of did he learn it from earlier sources? He does not say. (In fact, contrary to the rules of scholarship today, by which one must scrupulously indicate sources and give credit where it is due, it seems to have been the custom in those days not to indicate sources. For example, Descartes in Book III of his Géométrie explains how any cubic or quartic equation can be solved by taking square roots, cube roots, or trisecting an angle: this result is Prop. XXV of Viète's Supplementum Geometriae , but Descartes states it without giving any credit to Viète. )\hfill \break

This same trisection appears in the Mathematical Collection of Pappus [....] and in the book of lemmas attributed to Archimedes [Archimedes, 1792, Prop VIII p. 358]. Now the latter is known only through an Arabic source, and was not published until 1657 [,,,], so we may assume that Viète was unaware of it. However, Pappus' Mathematical Collection was circulating in manuscript in Europe well before the time of Viète, and was published in a Latin translation by Commandino in 1588 [see intro. to Pappus, 1986], so we may presume that Viète was familiar with the works of Pappus, and learned this construction from his book. (Commentators have also assumed that  borrowed the terms zetetics and poristics from Pappus, Book VII-see [Viète , 1983, p. 11, note]).\hfill \break

Proposition XVI \cite{Shooten1646} gives in geometrical form the solution of a cubic equation using the trisection of an angle. The statement is as follows:\hfill \break

If two isosceles triangles have their sides equal, and the angle at the base of the second is triple the angle at the base of the first, then the cube of the base of the first minus triple the solid made by the base of the first and the square of the common side, is equal to the solid on the base of the second and the same square of the common side.\hfill \break

Proposition XVI \cite{Shooten1646} gives in geometrical form the solution of a cubic equation using the trisection of an angle. The statement is as follows:\hfill \break

If two isosceles triangles have their sides equal, and the angle at the base of the second is triple the angle at the base of the first, then the cube of the base of the first minus triple the solid made by the base of the first and the square of the common side, is equal to the solid on the base of the second and the same square of the common side.\hfill \break
cm Let the first triangle be ABC with $AB=BC$, and the second $CDE$ with $CD=DE=BC$, and assume that the angle $DCE$ at the base of the second is three times the angle at A. Then he says
$$AC^3-3AC\cdot BC^2 = CE\cdot BC^2$$

Or, if we take the radius $BC=1$, let $AC=x$ and $CE=b$, the equation is

$$x^3-3x = b$$

Because the angle at A is one third of the angle $DCE$, the points $A,B,D$ are collinear and we can use the previous diagram. Viète's proof is expressed all in words. Translated into symbols it goes like this: drop perpendiculars from $B$ and $D$ to $ACE$, and let $y=BD$. By Euclid (III.36) we have

$$1\cdot (1+y) = (x-1)(x+1),$$

using the point $A$ and the two secants. Because the perpendiculars from $B$ and $D$ to $ACE$ are parallel, we have a proportion (VI.2)
$$1:y=\frac{1}{2}x : \frac{1}{2}(x+b)$$

Cross multiplying,

$$x+b=xy.$$

Now, solving the first equation for y and substituting, we find

$$x^3-3x = b$$

as required.\hfill \break

From the proposition of Viète, we see that to solve the equation  $x^3-3x = b$ , we make an isosceles triangle with sides 1 and base $b$, trisect the angle, and then the base of the new isosceles triangle will be $x$. By suitable linear substitutions, this gives a recipe for solving any cubic equation with three real roots, the so-called ``casus irreducibilis''. This is the case in which Cardano's formulas require cube roots of imaginary quantities, which were not understood at the time. As far as I know, this result is original to Viète, and this is the place where he first explains it.\hfill \break

Viète gives the general rule for this case of the cubic equation in De Recognitione Equationum \cite{Shooten1646}, but without proof, referring to his book Theoremata ad sectiones angulares  \cite{Shooten1646} for explanation.\hfill \break

Now we come to the construction of the regular heptagon. Viète's solution is directly analogous to Euclid's construction of the pentagon. Euclid first divides a line $AB$ in extreme and mean ratio: that is he finds a point C for which $AB:AC=AC:BC$. (II.11). From this he constructs an isosceles triangle whose base angles are twice the angle at the vertex (IV.10). Then from this triangle he constructs the pentagon (IV.11).\hfill \break

Viète first extends the diameter $BAC$ of a circle by a point I so that  $IB:IA=AB^2:IC^2$ (Prop. XIX). Form this he constructs an isosceles triangle whose base angles are three times the vertex angle (Prop XXIII), with which he can construct the heptagon (Prop XXIV).\hfill \break

Here is the construction of I. Let the diameter of the circle be $BAC$, and call the radius 1. Take $D$ with $AD=(\frac{1}{3})$. Let E be a vertex of the hexagon inscribed in the circle, so that $CE$ is $\frac{1}{6}$ of the circle. Draw $DE$. Make $AF$ parallel to $DE$. Draw $FGH$ so that $GH=AB$ (neusis, trisecting the angle $FAC$). Draw $EI$ parallel to $FG$. Then $I$ is the required point. For the proof, draw $AH$. Make $DK$ parallel to $AH$. Take $L$= midpoint of $DC$, and $M$= midpoint of $DL$. Then the triangles $IKD$ and $DEL$ satisfy the hypotheses of the previous result. So, if we let $x=ID$, $r=DE$ and note $DL=\frac{1}{3}$ ,then x will satisfy the equation
$$x^3-3xr^2=\frac{1}{3}r^2.$$

Now apply (I.47) to the triangle $DEM$. The base is  $\frac{1}6{}$ and the height is  $\frac{1}{2}\sqrt{3}$ , since $ACE$ is an equilateral triangle. Thus we find  $r^2=\frac{7}{9}$ , so $x$ satisfies

$$x^3-\frac{7}{3}x=\frac{7}{27}.$$

To verify the relation  $IB\cdot IC^2=IA\cdot AB^2$ , just note  $IB=x-\frac{4}{3}$,$IC=x+\frac{2}{3}$,$IA=x-\frac{1}{3}$ and $AB=1$. Substituting gives the equation above, so $I$ has the required property.\hfill \break

Now, changing notation, let $BAC$ be the diameter of a circle, let $D$ be the point $I$ as above, and draw $DE=AB=1$. Let $DE$ cut the circle again at $F$, and draw $AF$. Then $I$ claim $BE$ is the side of the regular heptagon inscribed in the circle. Draw $EC$. By (III.36),

$$DF\cdot DE=DB\cdot DC,$$

but also by construction

$$DB\cdot DC^2=DA\cdot AB^2$$

Since $DE=1=AB$, we find  $DF\cdot DC=DA\cdot DE$ . As a proportion, this says $DF:DE=DA:DC$. It follows that EC is parallel to AF (VI.2). Hence the angle at $C$ is equal to the angle at $BAF$. Denote it by 1. On the other hand, the angle at C subtends the arc $BE$, so the angle at the center $BAE$ is 2. So angle $EAF$ is equal to 1, and we find $AF$ is the angle bisector of $BAE$. Since $DAE$ is isosceles, the angle at $D$ is 2, and the exterior angle of the triangle $DAF$ at F is 3. But $AEF$ is isosceles, so the angle at $AEF$ is also 3. Now we have an isosceles triangle $AEF$ whose base angles are three times the vertex angle. It follows the vertex angle 1 is one seventh of two right angles. The angle $BAE$ is twice that, so $BE$ is the side of the regular heptagon, as required.\hfill \break

Note. In this account I have shortened Viète's exposition considerably by using modern algebraic notation, and in a few places I have done the calculation slightly differently, but in no way have I altered the basic mathematical content of the proof.\hfill \break

History of constructions of the heptagon. There is an Arabic text, discovered only in 1927 containing a construction of the heptagon attributed to Archimedes, which depends on the following: Given a square $ABCD$ with its diagonal $AD$, find a line $CE$ so that the triangles $CFD$ and $BGE$ have the same area  \cite{Waerden1954}, \cite{Knor1986} .This construction could be effected with an intersection of two conics, but it is not clear how to do it with a neusis.\hfill \break

Campanus, in his commentary to Euclid, Book IV, mentions that you could construct a heptagon if you had an isosceles triangle whose base angles are three times the vertex angle. \cite{Euklidis} Otherwise, this construction of the heptagon seems to be totally original to Viète.\hfill \break

From a modern point of view, to construct a heptagon one must solve the cyclotomic equation $x^7-1=0$  . If $\cos{\frac{2\pi}{7}}=\zeta+\zeta^{-1}$  is a root, then  is a root of the cubic equation
$$y^3+y^2-2y-1=0$$

One knows (from Viète!) that such a cubic can be solved by trisecting a certain angle, and modern constructions have been given following this idea by Plemelj, Bieberbach, Gleason and Hartshorne \cite{Hartshorne1997}. But none of these authors seems to have been aware of Viète's construction, almost 400 years earlier, based on the same ideas, but without the benefit of the modern theory of the cyclotomic equation and complex numbers.

\section{Adriaan van Roomen's problem}
This famous controversy is told by Tallemant des Réaux in these terms (46th story from the first volume of Les Historiettes. Mémoires pour servir à l’histoire du XVIIe siècle):\\

"In the times of Henri the fourth, a Dutchman called Adrianus Romanus, a learned mathematician, but not so good as he believed, published a treatise in which he proposed a question to all the mathematicians of Europe, but did not ask any Frenchman. Shortly after, a state ambassador came to the King at Fontainebleau. The King took pleasure in showing him all the sights, and he said people there were excellent in every profession in his kingdom. 'But, Sire,' said the ambassador, 'you have no mathematician, according to Adrianus Romanus, who didn't mention any in his catalog.' 'Yes, we have,' said the King. 'I have an excellent man. Go and seek Monsieur Viette,' he ordered. Vieta, who was at Fontainebleau, came at once. The ambassador sent for the book from Adrianus Romanus and showed the proposal to Vieta, who had arrived in the gallery, and before the King came out, he had already written two solutions with a pencil. By the evening he had sent many other solutions to the ambassador."\\

This suggests that the Adrien van Roomen problem is an equation of 45°, which Viète recognized immediately as a chord of an arc of 8° (${\frac  {2\pi }{45}}$ radians). It was then easy to determine the following 22 positive alternatives, the only valid ones at the time.\\

When, in 1595, Viète published his response to the problem set by Adriaan van Roomen, he proposed finding the resolution of the old problem of Apollonius, namely to find a circle tangent to three given circles. Van Roomen proposed a solution using a hyperbola, with which Viète did not agree, as he was hoping for a solution using Euclidean tools.\\

Viète published his own solution in 1600 in his work Apollonius Gallus. In this paper, Viète made use of the center of similitude of two circles.\cite{Cantor1911} His friend De Thou said that Adriaan van Roomen immediately left the University of Würzburg, saddled his horse, and went to Fontenay-le-Comte, where Viète lived. According to De Thou, he stayed a month with him and learned the methods of the new algebra. The two men became friends and Viète paid all van Roomen's expenses before his return to Würzburg.\\

This resolution had an almost immediate impact in Europe and Viète earned the admiration of many mathematicians over the centuries. Viète did not deal with cases (circles together, these tangents, etc.), but recognized that the number of solutions depends on the relative position of the three circles and outlined the ten resulting situations. Descartes completed (in 1643) the theorem of the three circles of Apollonius, leading to a quadratic equation in 87 terms, each of which is a product of six factors (which, with this method, makes the actual construction humanly impossible)\cite{Henk2003}
\section{Publications}
\subsection{Chronological list}
\begin{itemize}
    \item 

Between 1564 and 1568, Viète prepared for his student, Catherine de Parthenay, some textbooks of astronomy and trigonometry and a treatise that was never published: Harmonicon coeleste.

\item In 1571, Francisci Vietaei Universalium inspectionum ad Canonem mathematicum liber singularis (a book of trigonometry, often abbreviated Canonem mathematicum), which he published at his own expense and with great printing difficulties. This text contains many formulas on the sine and cosine and is unusual in using decimal numbers. The trigonometric tables here exceeded those of Regiomontanus (Triangulate Omnimodis, 1533) and Rheticus (1543, annexed to De revolutionibus of Copernicus).

\item In 1589, Deschiffrement d'une lettre escripte par le Commandeur Moreo au Roy d'Espaigne son maître.

\item In 1590, Deschiffrement description of a letter by the Commander Moreo at Roy Espaigne of his master, Tours: Mettayer.

\item In 1591:
\begin{itemize}
    \item 

In artem analyticem isagoge (Introduction to the art of analysis), also known as Algebra Nova (New Algebra) Tours: Mettayer, in 9 folio; the first edition of the Isagoge.
\item Zeteticorum libri quinque. Tours: Mettayer, in 24 folio; which are the five books of Zetetics, a collection of problems from Diophantus solved using the analytical art.
\end{itemize}

\item Between 1591 and 1593, Effectionum geometricarum canonica recensio. Tours: Mettayer, in 7 folio.

\item In 1593:
\begin{itemize}
    \item 

Vietae Supplementum geometriae. Tours: Francisci, in 21 folio.
\item Francisci Vietae Variorum de rebus responsorum mathematics liber VIII. Tours: Mettaye, in 49 folio; about the challenges of Scaliger.
\item Variorum de rebus mathematicis responsorum liber VIII; the "Eighth Book of Varied Responses" in which he talks about the problems of the trisection of the angle (which he acknowledges that it is bound to an equation of third degree) of squaring the circle, building the regular heptagon, etc.
\end{itemize}

\item In 1594, Munimen adversus nova cyclometrica. Paris: Mettayer, in quarto, 8 folio; again, a response against Scaliger.

\item In 1595, Ad problema quod omnibus mathematicis totius orbis construendum proposuit Adrianus Romanus, Francisci Vietae responsum. Paris: Mettayer, in quarto, 16 folio; about the Adriaan van Roomen problem.

\item In 1600:
\begin{itemize}
    \item 

De numerosa potestatum ad exegesim resolutione. Paris: Le Clerc, in 36 folio; work that provided the means for extracting roots and solutions of equations of degree at most 6.
\item Francisci Vietae Apollonius Gallus. Paris: Le Clerc, in quarto, 13 folio; where he referred to himself as the French Apollonius.
\end{itemize}

\item Between 1600 and 1602:
\begin{itemize}
    \item

Fontenaeensis libellorum supplicum in Regia magistri relatio Kalendarii vere Gregoriani ad ecclesiasticos doctores exhibita Pontifici Maximi Clementi VIII. Paris: Mettayer, in quarto, 40 folio.
\item Francisci Vietae adversus Christophorum Clavium expostulatio. Paris: Mettayer, in quarto, 8 folio; his theses against Clavius.
\end{itemize}

\end{itemize}

\subsection{Posthumous publications}
\begin{itemize}
    \item 
1612:
\begin{itemize}
    \item 

Supplementum Apollonii Galli edited by Marin Ghetaldi.
\item Supplementum Apollonii Redivivi sive analysis problematis bactenus desiderati ad Apollonii Pergaei doctrinam a Marino Ghetaldo Patritio Regusino hujusque non ita pridem institutam edited by Alexander Anderson.
\end{itemize}
\item1615:
\begin{itemize}
    \item 

Ad Angularum Sectionem Analytica Theoremata F. Vieta primum excogitata at absque ulla demonstratione ad nos transmissa, iam tandem demonstrationibus confirmata edited by Alexander Anderson.
\item Pro Zetetico Apolloniani problematis a se jam pridem edito in supplemento Apollonii Redivivi Zetetico Apolloniani problematis a se jam pridem edito; in qua ad ea quae obiter inibi perstrinxit Ghetaldus respondetur edited by Alexander Anderson
\item Francisci Vietae Fontenaeensis, De aequationum — recognitione et emendatione tractatus duo per Alexandrum Andersonum edited by Alexander Anderson
\end{itemize}
\item 1617: Animadversionis in Franciscum Vietam, a Clemente Cyriaco nuper editae brevis diakrisis edited by Alexander Anderson
\item 1619: Exercitationum Mathematicarum Decas Prima edited by Alexander Anderson
\item1631: In artem analyticem isagoge. Eiusdem ad logisticem speciosam notae priores, nunc primum in lucem editae. Paris: Baudry, in 12 folio; the second edition of the Isagoge, including the posthumously published Ad logisticem speciosam notae priores.
\end{itemize}

\section{Reception and influence}
During the ascendancy of the Catholic League, Viète's secretary was Nathaniel Tarporley, perhaps one of the more interesting and enigmatic mathematicians of 16th-century England. When he returned to London, Tarporley became one of the trusted friends of Thomas Harriot.\\

Apart from Catherine de Parthenay, Viète's other notable students were: French mathematician Jacques Aleaume, from Orleans, Marino Ghetaldi of Ragusa, Jean de Beaugrand and the Scottish mathematician Alexander Anderson. They illustrated his theories by publishing his works and continuing his methods. At his death, his heirs gave his manuscripts to Peter Aleaume.\cite{Auguste2008} We give here the most important posthumous editions:
\begin{itemize}
    \item 

In 1612: Supplementum Apollonii Galli of Marino Ghetaldi.
\item From 1615 to 1619: Animadversionis in Franciscum Vietam, Clemente a Cyriaco nuper by Alexander Anderson
\item Francisci Vietae Fontenaeensis ab aequationum recognitione et emendatione Tractatus duo Alexandrum per Andersonum. Paris, Laquehay, 1615, in 4, 135 p. The death of Alexander Anderson unfortunately halted the publication.
\item In 1630, an Introduction en l'art analytic ou nouvelle algèbre ('Introduction to the analytic art or modern algebra),\cite{Witmer1983} translated into French and commentary by mathematician J. L. Sieur de Vaulezard. Paris, Jacquin.
\item The Five Books of François Viette's Zetetic (Les cinq livres des zététiques de François Viette), put into French, and commented increased by mathematician J. L. Sieur de Vaulezard. Paris, Jacquin, p. 219.
\end{itemize}
The same year, there appeared an Isagoge by Antoine Vasset (a pseudonym of Claude Hardy), and the following year, a translation into Latin of Beaugrand, which Descartes would have received.\\

In 1648, the corpus of mathematical works printed by Frans van Schooten, professor at Leiden University (Elzevirs presses). He was assisted by Jacques Golius and Mersenne.\\

The English mathematicians Thomas Harriot and Isaac Newton, and the Dutch physicist Willebrord Snellius, the French mathematicians Pierre de Fermat and Blaise Pascal all used Viète's symbolism.\\

About 1770, the Italian mathematician Targioni Tozzetti, found in Florence Viète's Harmonicon coeleste. Viète had written in it: Describat Planeta Ellipsim ad motum anomaliae ad Terram. (That shows he adopted Copernicus' system and understood before Kepler the elliptic form of the orbits of planets.) \cite{Swerdlow1975} \\

In 1841, the French mathematician Michel Chasles was one of the first to reevaluate his role in the development of modern algebra.\\

In 1847, a letter from François Arago, perpetual secretary of the Academy of Sciences (Paris), announced his intention to write a biography of François Viète.\\

Between 1880 and 1890, the polytechnician Fréderic Ritter, based in Fontenay-le-Comte, was the first translator of the works of François Viète and his first contemporary biographer with Benjamin Fillon.

\section{Descartes' views on Viète}
Thirty-four years after the death of Viète, the philosopher René Descartes published his method and a book of geometry that changed the landscape of algebra and built on Viète's work, applying it to the geometry by removing its requirements of homogeneity. Descartes, accused by Jean Baptiste Chauveau, a former classmate of La Flèche, explained in a letter to Mersenne (1639 February) that he never read those works.[23] Descartes accepted the Viète's view of mathematics for which the study shall stress the self-evidence of the results that Descartes implemented translating the symbolic algebra in geometric reasoning.\cite{Bullynk2018} The locution mathesis universalis was derived from van Roomen's works.\cite{Bullynk2018} \hfill \break

"I have no knowledge of this surveyor and I wonder what he said, that we studied Viète's work together in Paris, because it is a book which I cannot remember having seen the cover, while I was in France."\hfill \break

Elsewhere, Descartes said that Viète's notations were confusing and used unnecessary geometric justifications. In some letters, he showed he understands the program of the Artem Analyticem Isagoge; in others, he shamelessly caricatured Viète's proposals. One of his biographers, Charles Adam,\cite{Adam} noted this contradiction:\hfill \break

"These words are surprising, by the way, for he (Descartes) had just said a few lines earlier that he had tried to put in his geometry only what he believed "was known neither by Viète nor by anyone else".\hfill \break

So he was informed of what Viète knew; and he must have read his works previously."
Current research has not shown the extent of the direct influence of the works of Viète on Descartes. This influence could have been formed through the works of Adriaan van Roomen or Jacques Aleaume at the Hague, or through the book by Jean de Beaugrand.\cite{Sasaki2013} \hfill \break

In his letters to Mersenne, Descartes consciously minimized the originality and depth of the work of his predecessors. "I began," he says, "where Viète finished". His views emerged in the 17th century and mathematicians won a clear algebraic language without the requirements of homogeneity. Many contemporary studies have restored the work of Parthenay's mathematician, showing he had the double merit of introducing the first elements of literal calculation and building a first axiomatic for algebra.\cite{Hairer2008} \hfill \break

Although Viète was not the first to propose notation of unknown quantities by letters - Jordanus Nemorarius had done this in the past - we can reasonably estimate that it would be simplistic to summarize his innovations for that discovery and place him at the junction of algebraic transformations made during the late sixteenth – early 17th century.

\bibliographystyle{unsrtnat}
\bibliography{references}  

\begin{thebibliography}{24}
\providecommand{\natexlab}[1]{#1}
\providecommand{\url}[1]{\texttt{#1}}
\expandafter\ifx\csname urlstyle\endcsname\relax
  \providecommand{\doi}[1]{doi: #1}\else
  \providecommand{\doi}{doi: \begingroup \urlstyle{rm}\Url}\fi

\bibitem[Chisholm(1911)]{Cantor1911}
Hugh Chisholm.
\newblock \emph{The Encyclop{\ae}dia Britannica: Vetch-Zymotic Diseases},
  volume~28.
\newblock At the University Press, 1911.

\bibitem[Britannica(2016)]{Brit2016}
Encyclopedia Britannica, 2016.
\newblock entry for Mathematics in the 17th and 18th centuries.

\bibitem[Kline(1990)]{Kline1972}
Morris Kline.
\newblock \emph{Mathematical Thought from Ancient to Modern Times: Volume 2},
  volume~2.
\newblock Oxford university press, 1990.

\bibitem[Vi{\`e}te(1970{\natexlab{a}})]{Vonwort1970}
Fran{\c{c}}ois Vi{\`e}te.
\newblock Opera mathematica, recognita: Francisci a schooten, vorwort und
  register von joseph e.
\newblock \emph{Hofmann, New York: G. Olms}, 1970{\natexlab{a}}.

\bibitem[Vi{\`e}te and Witmer(2006)]{Witmer1983}
Fran{\c{c}}ois Vi{\`e}te and T~Richard Witmer.
\newblock \emph{The analytic art}.
\newblock Courier Corporation, 2006.

\bibitem[von Camerer et~al.(1795)]{Camerer1795}
Johann~Wilhelm von Camerer et~al.
\newblock \emph{De tactionibus quae supersunt...}
\newblock Ettinger, 1795.

\bibitem[Vi{\`e}te(1970{\natexlab{b}})]{Shooten1646}
Fran{\c{c}}ois Vi{\`e}te.
\newblock Opera mathematica" published by franciscus van schooten with
  commentaries, leiden, 1646.
\newblock \emph{Reprint by Georg Olms Verlag, Hildesheim, New York},
  1970{\natexlab{b}}.

\bibitem[Descartes and {\v{Z}}i{\v{z}}ek(2013)]{Cogito2013}
Rene Descartes and Slavoj {\v{Z}}i{\v{z}}ek.
\newblock \emph{Cogito ergo sum}.
\newblock AlboVersorio, 2013.

\bibitem[Collins()]{Connor}
Lincoln Collins.
\newblock Biography paper: Francois viete.

\bibitem[Charbonneau(2003)]{Brit}
Louis Charbonneau.
\newblock Il ya 400 ans mourait sieur fran{\c{c}}ois vi{\`e}te, seigneur de la
  bigoti{\`e}re.
\newblock \emph{Bulletin AMQ}, 2003.

\bibitem[Bos(2001)]{Bos}
Henk~JM Bos.
\newblock \emph{Redefining geometrical exactness: Descartes' transformation of
  the early modern concept of construction}.
\newblock Springer Science \& Business Media, 2001.

\bibitem[Pycior and Seltman(1998)]{Pycior}
Helena~M Pycior and M~Seltman.
\newblock Reviews: Mathematics and logic-symbols, impossible numbers, and
  geometric entanglements: British algebra through the commentaries on newton's
  universal arithmetick.
\newblock \emph{Annals of Science}, 55\penalty0 (4), 1998.

\bibitem[Roegel(2021)]{Canon1579}
Denis Roegel.
\newblock \emph{A reconstruction of Vi{\`e}te's Canon Mathematicus (1579)}.
\newblock PhD thesis, LORIA (Universit{\'e} de Lorraine, CNRS, INRIA), 2021.

\bibitem[Van~der Waerden(1954)]{Waerden1954}
BL~Van~der Waerden.
\newblock Science awakening, p. noordhoff ltd-groningen holland.
\newblock \emph{English Translation by Arnold Dresden}, 1954.

\bibitem[Knorr(1986)]{Knor1986}
Wilbur~Richard Knorr.
\newblock The ancient tradition of greek geometric problems.
\newblock \emph{Basel, Boston, and Berlin: Birkh{\"a}user Verlag}, 1986.

\bibitem[Euclides et~al.()Euclides, Veneto, Briconneto, and
  Hypsicles]{Euklidis}
Euclides, Bartholomaeo~Zamberto Veneto, Francisco Briconneto, and Hypsicles.
\newblock \emph{Euclidis Megarensis Geometricorum elementorum libri XV}.
\newblock Henricus Stephanus.

\bibitem[Hartshorne(1997)]{Hartshorne1997}
Robin Hartshorne.
\newblock Companion to euclid: A course of geometry.
\newblock \emph{Based on Euclid’s Elements and Its Modern Descendants,
  Berkeley: American Mathematical Society}, 1997.

\bibitem[Bos(2003)]{Henk2003}
Henk~JM Bos.
\newblock Descartes, elisabeth and apollonius’ problem.
\newblock \emph{The Correspondence of Ren{\'e} Descartes}, 1643:\penalty0
  202--212, 2003.

\bibitem[Thou(2008)]{Auguste2008}
Jacques-Auguste~de Thou.
\newblock La vie de jacques-auguste de thou. i. aug. thuani vita.
\newblock 2008.

\bibitem[Swerdlow(1975)]{Swerdlow1975}
Noel~M Swerdlow.
\newblock The planetary theory of fran{\c{c}}ois vi{\`e}te 1. the fundamental
  planetary models.
\newblock \emph{Journal for the History of Astronomy}, 6\penalty0 (3):\penalty0
  185--208, 1975.

\bibitem[Bullynck(2018)]{Bullynk2018}
Maarten Bullynck.
\newblock The'everyday'in mathematics on the usability of mathematical
  practices for doing history.
\newblock 2018.

\bibitem[Adam(1910)]{Adam}
Charles Adam.
\newblock \emph{Vie \& oeuvres de Descartees: {\'e}tude historique}, volume~12.
\newblock L{\'e}opold Cerf, 1910.

\bibitem[Sasaki(2013)]{Sasaki2013}
Chikara Sasaki.
\newblock \emph{Descartes’s mathematical thought}, volume 237.
\newblock Springer Science \& Business Media, 2013.

\bibitem[Hairer and Wanner(2008)]{Hairer2008}
E~Hairer and G~Wanner.
\newblock Introduction to analysis of the infinite.
\newblock In \emph{Analysis by Its History}, pages 1--79. Springer, 2008.

\end{thebibliography}






\end{document}